\documentclass[letterpaper]{article}
\usepackage[intlimits]{amsmath}
\usepackage{amsfonts, amssymb, amsthm}
\newlength{\hchng}
\newlength{\vchng}
\newtheorem{thm}{Theorem}[section]

\newtheorem{cor}[thm]{Corollary}
\newtheorem{lemma}[thm]{Lemma}

\newtheorem{preremark}[thm]{Remark}
\newenvironment{remark}{\begin{preremark}\rm}{\medskip \end{preremark}}
\numberwithin{equation}{section}

\newcommand{\abs}[1]{\left\vert#1\right\vert}

\newcommand{\R}{\mathbb R}
\newcommand{\eps}{\varepsilon}

\newcommand{\grad} {\nabla}
\newcommand{\lap} {\triangle}
\newcommand{\bdary} {\partial}
\newcommand{\dx} {\; \mathrm{d} x}
\newcommand{\dd} {\; \mathrm{d}}

\DeclareMathOperator{\dv}{div}

\newcommand{\Ly}[1]{\lap_x #1 + \frac{a}{y} #1_y + #1_{yy}}
\newcommand{\Dy}[1]{\dv \left(y^a \grad #1\right)}
\newcommand{\Lz}[1]{\lap_x #1 + z^\alpha #1_{zz}}

\newcommand{\party}[1] {#1_y}
\newcommand{\partyy}[1] {#1_{yy}}

\newcommand{\partx}[1] {#1_x}
\newcommand{\partxx}[1] {#1_{xx}}

\newcommand{\partz}[1] {#1_z}
\newcommand{\partzz}[1] {#1_{zz}}

\newcommand{\ydd} {\; y^a \mathrm{d}}

\title{An extension problem related to the fractional Laplacian}
\author{Luis Caffarelli and Luis Silvestre}
\begin{document}

\maketitle

\begin{abstract}
The operator square root of the Laplacian $(-\lap)^{1/2}$ can be obtained from the harmonic extension problem to the upper half space as the operator that maps the Dirichlet boundary condition to the Neumann condition. In this paper we obtain similar characterizations for general fractional powers of the Laplacian and other integro-differential operators. From those characterizations we derive some properties of these integro-differential equations from purely local arguments in the extension problems.
\end{abstract}

\section{Introduction}

Let us suppose we have a smooth bounded function $f : \R^n \to \R$ and we solve the extension problem
\begin{align}
u(x,0) &= f(x) && \text{for } x \in \R^n \label{eq:dc}\\
\lap u(x,y) &= 0 && \text{for } x \in \R^n \text{ and } y > 0 \label{eq:lap}
\end{align}
to obtain a smooth bounded function $u$. It is well known that $-u_y (x,0) = (-\lap)^{1/2} f(x)$, and therefore we can realize $(-\lap)^{1/2}$ as the operator $T : f \mapsto -\party u (x,0)$ in the above extension problem.

This is easy to show by applying $T$ twice. When we place $-u_y(x,0)$ instead of $f$ as the Dirichlet condition in \eqref{eq:dc}, we obtain $-u_y(x,y)$ instead of $u$ as the solution of \eqref{eq:dc}-\eqref{eq:lap}. Then $T(T(f))(x) = T(-u_y(x,0))(x) = u_{yy}(x,0) = - \lap_{x} f(x)$. To show $T = (-\lap)^{1/2}$ it is only left to check that $T$ is indeed a positive operator, which follows by a simple integration by parts argument.

In this paper we will generalize this situation to a similar extension problem for each fractional power of the Laplacian. We will construct any fractional Laplacian from an extension problem to the upper half space for a specific elliptic partial differential equation. This allows us (in forthcoming work) to treat non linear variational problems, involving fractional Laplacians, with standard local perturbation methods from the calculus of variations. Examples of such methods are the Almgren Frequency formula, and the Boundary Harnack inequality presented below. The partial differential equation that we obtain will turn out to be degenerate for any power of the Laplacian other than $(-\lap)^{1/2}$, however they belong to a more general class of equations that shares many of the essential properties of uniformly elliptic equations (as in \cite{FKS}, \cite{FJK}, \cite{FKJ} for the divergence case, or \cite{CG} for the non-divergence case). %We will study the properties of this degenerate elliptic equation. Then we will derive a few properties of the fractional Laplacian like Harnack and boundary Harnack inequalities from older results in degenerate elliptic equations \cite{FKS}, \cite{FJK}, \cite{FKJ}, \cite{CG}. Our ultimate goal is to use these properties for studying free boundary problems for the fractional Laplacian, in particular to extend the result of \cite{S}, in future papers.

The fractional Laplacian of a function $f : \R^n \to \R$ is expressed by the formula
\begin{equation}
(-\lap)^s f(x) = C_{n,s} \int_{\R^n} \frac{f(x) - f(\xi)}{\abs{x-\xi}^{n+2s}} \dd \xi
\end{equation}
where the parameter $s$ is a real number between $0$ and $1$, and $C_{n,s}$ is some normalization constant.

It can also be defined as a pseudo-differential operator
\[ \widehat {(-\lap)^s f} (\xi) = \abs{\xi}^{2s} \widehat f(\xi) \]

The fractional Laplacian can be defined in a distributional sense for functions that are not differentiable as long as $\hat f$ is not too singular at the origin or, in terms of the $x$ variable, as long as
\[ \int_{\R^n} \frac{\abs{f(x)}}{(1+\abs{x})^{n+2s}} \dd x < + \infty \]

We will relate the fractional Laplacian with solutions of the following extension problem. For a function $f : \R^n \to \R$, we consider the extension $u : \R^n \times [0,\infty) \to \R$ that satisfies the equation
\begin{align}
u(x,0) &= f(x) \label{eq:dirichletboundary}\\
\Ly u &= 0 \label{eq:withy}
\end{align}

The equation \eqref{eq:withy} can also be written as
\begin{equation}
\Dy u = 0
\end{equation}

Which is clearly the Euler-Lagrange equation for the functional
\begin{equation} \label{eq:functional}
 J(u) = \int_{y>0} \abs{ \grad u }^2 \ydd X 
\end{equation}

We will show that
\[ C (-\lap)^s f = \lim_{y \to 0^+} -y^a u_y = \frac{1}{1-a} \lim_{y \to 0} \frac{u(x,y) - u(x,0)}{y^{1-a}}\]
for $s = \frac{1-a}{2}$ and some constant $C$ depending on $n$ and $s$. Which reduces to the regular normal derivative in the case $a=0$ (as in \eqref{eq:dc}-\eqref{eq:lap}).

If we make the change of variables $z = \left( \frac{y}{1-a} \right)^{1-a}$ in \eqref{eq:withy}, we obtain a nondivergence form equation
\begin{equation}\label{eq:withz}
\Lz u = 0
\end{equation}
for $\alpha = \frac{-2a}{1-a}$. Moreover, $y^a u_y = u_z$. Thus, we claim that the following equality holds up to a multiplicative constant
\[ (-\lap)^s f(x) = -\lim_{y \to 0^+} y^a u_y(x,y) = -u_z(x,0) \]

It seems convenient to keep this notation. The variable $x$ is always in $\R^n$. The variables $y$ and $z$ are nonnegative real numbers that satisfy the relation $z = \left( \frac{y}{1-a} \right)^{1-a}$. The function $u$ is the extension of $f$ to the upper half space and takes values in $\R^n \times [0,\infty)$. We will use $u$ and its version after the change of variables $y \mapsto z$ indistinctly, and we will call the variable either $y$ or $z$ to point out the difference. Whenever we refer to a point in $\R^{n+1}$ we will use capital letters (like $X$).

\section{Properties of the PDEs}
In this section we will study the basic properties of the equations \eqref{eq:withy} and \eqref{eq:withz}. We will develop explicit Poisson formulas among other properties.

\subsection{Harmonic functions in $n+1+a$ dimensions} \label{sec:rsymmetric}

The equation \eqref{eq:withy} has a curious intuition behind it that will help us obtain several properties of its solutions.

For a nonnegative integer $a$, suppose $u(x,y) : \R^n \times \R^{1+a} \to \R$ is radially symmetric in the $y$ variable, meaning that if $\abs{y} = \abs{y'} = r$, then $u(x,y) = u(x,y')$. We can think of $u$ as a function of $x$ and $r$, and in these variables write an expresion for its Laplacian:
\[ \lap u =\lap_x u + \frac{a}{r} u_r + u_{rr} \]

We thus obtain an identical expression of equation \eqref{eq:withy} (with $r$ instead of $y$, but that's just naming). As far as the expresion is concerned, there is no need to consider only integer values of $a$. We can then think of the equation \eqref{eq:withy} as the harmonic extension problem of $f$ in $1+a$ dimensions more. Although it is impossible to think of a meaning for $\R^{n+1+a}$ when $a$ is not an integer, the solutions of \eqref{eq:withy} will satisfy many properties common to harmonic functions.

\subsection{Fundamental solution at the origin}
To obtain the \emph{fundamental solution} for \eqref{eq:withy} at the origin, we just have to consider the fundamental solution of the Laplacian in $n+1+a$ dimensions. If $n-1+a > 1$, it reads
\begin{equation}
 \Gamma(X) = C_{n+1+a} \frac{1}{\abs{X}^{n-1+a}}
\end{equation}
where the constant $C_{n+1+a}$ is given by the formula $C_k = \pi^{k/2} \Gamma( k/2 -1) / 4$.

It can be verified as a direct computation that $\Gamma$ is a solution of \eqref{eq:withy} when $y \neq 0$ and indeed $\lim_{y \to 0^+} y^a u_y = -C \delta_0$ for some constant $C$ as we will see later. Notice that $\Gamma(x,0) = \frac{C}{\abs{x}^{n-1+a}}$ is the fundamental solution of the fractional Laplacian $(-\lap)^{\frac{1-a}{2}}$ for some appropriate constant $C$ depending on $n$ and $a$ (Recall $X = (x,y)$).

Since the equation \eqref{eq:withz} can be derived from \eqref{eq:withy} by just a change of variables, we can also change variables to obtain a corresponding fundamental solution
\begin{equation}
 \tilde \Gamma(x,z) = C_{n+1+a} \frac{1}{\left( \abs{x}^2 + (1-a)^2 \abs{z}^{2/(1-a)} \right)^{\frac{n-1+a}{2}}}
\end{equation}
that solves \eqref{eq:withz} when $z \neq 0$ and $u_z(x,z) \to -\delta_0$ as $z \to 0$.

\subsection{Conjugate equation}
We have seen that if $u$ is a solution of \eqref{eq:withy}, then the function $w(x,y) := y^a u_y(x,y)$ seems to carry interesting information. Indeed, this function satisfies the conjugate equation
\[ \lap_x w - \frac{a}{y} w_y + w_{yy} = 0\]
that is nothing but \eqref{eq:withy} with $-a$ instead of $a$. This property of $w$ can be checked by a direct computation:
\begin{align*}
\lap_x w - \frac{a}{y} w_y + \partyy w &= y^a \left( \partial_y \lap_x u - \frac{a^2}{y^2} \party u -  \frac{a}{y}  \partyy u + \frac{a(a-1)}{y^2} \party u + 2 \frac{a}{y} \partyy u + u_{yyy} \right) \\
&= y^a \left( \partial_y \lap_x u - \frac{a}{y^2} \party u + \frac{a}{y} \partyy u + u_{yyy} \right) \\
&= y^a \partial_y \left( \lap_x u + \frac{a}{y} \party u + \partyy u \right) = 0
\end{align*}

The intuition behind the above fact is simple to explain when $n=1$ ($x\in \R$ and $y \in  [0,+\infty)$). The function $w$ turns out to be the stream function related to $u_x$ in the following sense: if we set $v = \partx u$, then clearly $v$ is a solution of \eqref{eq:withy}, since the equation is invariant under translations in $x$. Thus we have $\dv (y^a \grad v) = 0$. The vector field $(y^a \party v, -y^a \partx v)$ is then irrotational and there is a function $w$ such that $\grad w = (-y^a \party v, y^a \partx v)$. This function $w$ is the stream function of $v$, and it satisfies the equation $\dv (y^{-a} \grad w) = 0$. Now we check that $w = y^a \party u(x,y)$, since \[\grad (y^a \party u) = \left( y^a u_{xy} , y^a \left( \frac{a}{y} \party u + \partyy u \right) \right) = (y^a \party v, -y^a \partx v) \]

\subsection{Poisson formula}
We want to develop a Poisson formula $P$ to explicitly solve \eqref{eq:dirichletboundary}-\eqref{eq:withy}.
\[ u(X) = \int_{\R^n} P(x-\xi,y) f(\xi) \dd \xi \]

The Poisson kernel $P$ must be a solution to \eqref{eq:withy} where $y>0$ and $\lim_{y \to 0} P(x,y) = \delta_0$. from the previous sections, we see that the correct choice would be to take $P(x,y) = -y^{-a} \partial_y \Gamma_{-a}(x,y)$. Therefore
\begin{equation} \label{eq:poisson}
P(x,y) = %(n+a-1) C_{n-1-a} 
C_{n,a} \frac{y^{1-a}}{\left( \abs{x}^2 + \abs{y}^2 \right)^{\frac{n+1-a}{2}}}
\end{equation}

The fact that $P$ is a solution to \eqref{eq:withy} where $y>0$ is a consequence of the fact that $\Gamma_{-a}$ is to the conjugate equation. Of course it can also be checked by a direct computation. Moreover, $P(x,y) = y^{-n} P(x/y,1)$, so $P$ must converge to a multiple of the Dirac delta as $y \to 0$, so $P$ is indeed the Poisson kernel.

The corresponding Poisson kernel for the equation \eqref{eq:withz} can be obtained either by a change of variables from \eqref{eq:poisson} or by computing $-\tilde P_z$:
\begin{equation}
\tilde P(x,z) = %(n+a-1) (1-a)^{1/(1-a)} C_{n-1-a} 
C_{n,a} \frac{z}{\left(\abs{x}^2 + (1-a)^2 \abs{z}^{2/(1-a)} \right)^{\frac{n+1-a}{2}}}
\end{equation}

%Naturally, it is not hard to check that $P$ and $\tilde P$ are indeed the correct Poisson kernels by a straight forward computation.

\section{Relation with fractional Laplacian}
In this section we will see how the equations \eqref{eq:withy} or \eqref{eq:withz} relate to the operator $(-\lap)^s$. Namely, we will show that up to a constant factor
\begin{equation} \label{eq:identities}
 \lim_{y \to 0} y^a \party u(x,y) = \partz u(x,0) = -(-\lap)^s f(x) = \int_{\R^n} \frac{ f(\xi) - f(x) }{\abs{\xi-x}^{n+2s}} \dd \xi
\end{equation}

\subsection{Proof using the Poisson formula}
%We can compute $\lim_{y \to 0} (u(x,y)-u(x,0))/y^{1-a}$ using the Poisson formula
%\begin{align*}
%\lim_{y \to 0} \frac{u(x,y)-u(x,0)}{y^{1-a}} &= C_{n,a} \int_{\R^n} \frac{f(z)-f(x)}{\left( \abs{x-z}^2 + \abs{y}^2 \right)^{\frac{n+1-a}{2}}} \dd z \\
%&= C_{n,a} \; \mathrm{PV} \int_{\R^n} \frac{f(z)-f(x)}{ \abs{x-z}^{\frac{n+1-a}{2}}} \dd z\\
%&= -C(-\lap)^{\frac{1-a}{2}} f(x)
%\end{align*}
%as long as $f$ is regular enough.

We can compute $\partz u(x,0)$ using the Poisson formula.
\begin{align}
 \partz u(x,0) &= \lim_{z \to 0} \frac{u(x,z) - u(x,0)} {z} \\
 &= \lim_{z \to 0} \frac{1}{z} \int_{\R^n} \tilde P(x-\xi,z) (f(\xi) - f(x)) \dd \xi \\
 &= \lim_{z \to 0} \int_{\R^n} \frac{C}{\left(\abs{x-\xi}^2 + (1-a)^2 \abs{z}^{2/(1-a)} \right)^{\frac{n+1-a}{2}}} (f(\xi) - f(x)) \dd \xi \label{aux:withlim} \\
 &= C \mathrm{PV} \int_{\R^n} \frac{f(\xi) - f(x)}{\abs{x-\xi}^{n+1-a}} \dd \xi \label{aux:withoutlim}\\
 &= -C(-\lap)^{\frac{1-a}{2}} f(x)
\end{align}
where the limit in \eqref{aux:withlim} exists and equals \eqref{aux:withoutlim} as long as $f$ is regular enough.

Changing variables, this implies that $\lim_{y \to 0} -y^a u_y(x,y)$ also converges to a multiple of the fractional Laplacian. The fact that $y^a u_y(x,y)$ has a limit as $y \to 0$, immediately implies that $\lim_{y \to 0} \frac{u(x,y)-u(x,0)}{y^{1-a}}$ has the same limit.

\begin{remark}
The proof was done by computing $u_z$. If instead we had chosen to compute $y^a u_y$ or $\frac{u(x,y)-u(x,0)}{y^{1-a}}$ the complexity of the computation would be comparable.
\end{remark}

\subsection{Proof using Fourier transform} \label{sec:ft}
Alternatively, we can also prove \eqref{eq:identities} taking Fourier transform in $x$. 

One way to do this is proving that the corresponding energy functionals coincide. Namely
\begin{equation} \label{eq:energies}
 \int_{y>0} \abs{\grad u}^2 \ydd X = \int_{\R^n} \abs{\xi}^{2s} \abs{ \hat f(\xi) }^2 \dd \xi
\end{equation}

The equation \eqref{eq:withy} becomes
\[ -\abs{\xi}^2 \hat u(\xi,y) +\frac{a}{y} \hat u_y (\xi,y) + \hat u_{yy}(\xi,y) = 0 \]

We thus obtain an ordinary differential equation for each value of $\xi$. 

Suppose that $\phi : [0,\infty) \to \R$ is the minimizer of the functional
\[ J(\phi) := \int_{y>0} (\abs{\phi'}^2 + \abs{\phi}^2) \ydd y\]
for $\phi(0)=1$, then $\phi$ solves the following equation:
\begin{align*}
 -\hat \phi (y) +\frac{a}{y} \hat \phi_y(y) + \phi_{yy}(y) &= 0 \\
 \phi(0) &= 1 \\
 \lim_{y \to \infty} \phi(y) &=0
\end{align*}

By a simple scaling we can see that
\[ \hat u(\xi,y) = \hat f(\xi) \phi(\abs{\xi} y) \]
and thus, the energy of $u$ becomes
\begin{align*}
 \int_{y>0} \abs{\grad u}^2 \ydd X &= \int_{\R^n} \int_0^\infty  (\abs{\xi}^2 \abs{\hat u}^2 + \abs{\hat u_y}^2) y^a \dd y \dd \xi \\
 &= \int_{\R^n} \int_0^\infty \abs{f(\xi)}^2 \abs{\xi}^2 \left( \abs{\phi(\abs{\xi} y)}^2 + \abs{\phi'(\abs{\xi} y)}^2 \right) y^a \dd y \dd \xi \\
 &= \int_{\R^n} \abs{f(\xi)}^2 \abs{\xi}^{1-a} \int_0^\infty \left( \abs{\phi(\bar y)}^2 + \abs{\phi'(\bar y)}^2 \right) {\bar y}^a \dd \bar y \dd \xi \\
 &= \int_{\R^n} \abs{f(\xi)}^2 \abs{\xi}^{1-a} J(\phi) \dd \xi 
\end{align*}

Thus we conclude \eqref{eq:energies}. The corresponding Euler Lagrange equations for each energy must then coincide up to a constant factor. Therefore we obtain
\[ -\lim_{y \to 0} y^a \party u(x,y) = C (-\lap)^{\frac{1-a}{2}} f(x) \]
where the constant $C$ is given by $J(\phi)$.

We could also make another proof from \eqref{eq:withz}. After taking Fourier transform, the equation \eqref{eq:withz} becomes
\[ -\abs{\xi}^2 \hat u(\xi,z) + z^\alpha \hat u_{zz}(\xi,z) = 0 \]

We thus obtain an ordinary differential equation for each value of $\xi$.

Suppose that $\phi : [0,\infty) \to \R$ solves the following equation
\begin{equation} \label{eq:fwithz}
\begin{aligned}
 -z^\alpha \phi''(z) + \phi(z) &= 0 \\
 \phi(0) &= 1 \\
 \lim_{z \to \infty} \phi(z) &=0
\end{aligned}
\end{equation}

Then by a simple scaling we can see that
\[ \hat u(\xi,z) = \hat u(\xi,0) \phi(\abs{\xi}^{\frac{2}{2-\alpha}} z) \]

And then, 
\[ 
\begin{split}
\hat u_z(\xi,0) &= \hat u(\xi,0) \abs{\xi}^{\frac{2}{2-\alpha}} \phi'(0) \\
&= C_a \abs{\xi}^{1-a} \hat u(\xi,0) = C_a \abs{\xi}^{1-a} \hat f(\xi)
\end{split}
\]

To prove \eqref{eq:identities} we now only need to show that such function $\phi$ exists and it is differentiable at zero. We notice that for $\eps$ small enough then
\[ \overline \phi (z) := \min (1, z^{-\eps}) \]
is a supersolution, whereas for $A$ and $B$ large
\[ \underline{\phi}(z) := e^{-A z^{1/2} + B z^2} \]
is a supersolution. Thus, by Perron's method, we can find a $\phi$ in between that solves \eqref{eq:fwithz}. Moreover this $\phi$ is Lipschitz at $z=0$ and $0 \leq \phi \leq 1$ since it is trapped between $\underline{\phi}$ and $\overline{\phi}$. From the equation \eqref{eq:fwithz} we also see that then $0 \leq \partzz \phi \leq z^{-\alpha}$, thus $\phi \in C^{2-\alpha}$. Not only does this complete the proof of \eqref{eq:identities} but it also gives an interesting regularity estimate for the solution $u$ with respect to $z$. In the next section we will explore this issue further.

\section{Reflection extensions}
In order to apply interior Harnack estimates to our PDE \eqref{eq:withy} (or \eqref{eq:withz}), we must show that if the operator $(-\lap)^s f = 0$ in an open ball, then we can reflect the solution u and make it a solution of \eqref{eq:withy} (or \eqref{eq:withz}) across $y=0$ ($z=0$) in a suitable sense.

We first address the divergence case.
\begin{lemma} \label{lem:dref}
Suppose that $u : \R^n \times [0,\infty) \to \R$ is a solution to \eqref{eq:withy} such that for $|x| \leq r$
\begin{equation} \label{eq:zeron}
 \lim_{y \to 0} y^a \party u(x,y)=0
\end{equation}
then the extension to the whole space
\[ \tilde{u}(x,y) = \begin{cases}
u(x,y) & y \geq 0 \\
u(x,-y) & y<0
\end{cases} \]
is a solution to %\eqref{eq:withy} 
\[ \dv( |y|^a \grad u ) = 0 \]
in the week sense in the $(n+1)$ dimensional ball of radius $R$ ($\{ (x,y) : \abs{x}^2 + \abs{y}^2 \leq R^2 \}$).
\end{lemma}

\begin{proof}
Let $h \in C^\infty_0 (B_R)$ be a test function. We want to show that
\begin{equation} \label{eq:variational}
 \int_{B_R} \grad \tilde u \cdot \grad h \; |y|^a \dd X = 0
\end{equation}
Let $\eps > 0$, we separate a strip of width $\eps$ around $y=0$ in the domain of the integral above
\begin{align*}
 \int_{B_R} \grad \tilde u \cdot \grad h \; |y|^a \dd X &= \int_{B_R \setminus \{ \abs{y} < \eps \} } \grad \tilde u \cdot \grad h \; |y|^a \dd X + \int_{B_R \cap \{ \abs{y} < \eps \} } \grad \tilde u \cdot \grad h \; |y|^a \dd X \\
 &= \int_{B_R \setminus \{ \abs{y} < \eps \} } \dv \left( |y|^a h \grad \tilde u \right) \dd X + \int_{B_R \cap \{ \abs{y} < \eps \} } \grad \tilde u \cdot \grad h \; |y|^a \dd X \\
 &= \int_{B_R \cap \{ \abs{y} = \eps \} } h \tilde u_y (x,\eps) \; \eps^a \dd x + \int_{B_R \cap \{ \abs{y} < \eps \} } \grad \tilde u \cdot \grad h \; |y|^a \dd X
\end{align*}

When we let $\eps \to 0$, the second term of the right hand side above clearly goes to zero because $|y|^a \abs{\grad u}^2$ is locally integrable, and the first term converges to zero if $\eps^a \tilde u_y(x,\eps) \to 0$ as $\eps \to 0$.

Therefore, if $\eps^a \tilde u_y(x,\eps) \to 0$ as $\eps \to 0$, then $\dv(|y|^a \grad u)=0$ in the whole ball $B_R$ (accross $y=0$).
\end{proof}

We must clarify in what sense we take the limit in \eqref{eq:zeron}. In this particular case, since the limit vanishes, we could prove that actually $u$ is $C^\infty$ in $x$ and the limit holds in the uniform sense. In general it would be convenient to have a weak definition of \eqref{eq:zeron}. The equation \eqref{eq:variational} becomes the definition of \eqref{eq:zeron} in the weak sense. In any case the limit \eqref{eq:zeron} holds in the sense that, if we take any smooth test function $\varphi \in C^\infty_0 (B_1^*)$, \[\lim_{y \to 0} \int_{B_1^*} y^a u_y(x,y) \varphi(x) \dx = 0\]

For a general class of functions $g$, the problem
\begin{align*}
\dv (\abs{y}^a \grad u) &= 0 &&\text{in the weak sense in } B_R \\
u &= g && \text{on } \bdary B_R
\end{align*}
has a unique solution in the weighted Sobolev space $H^{1,2} (B_R, \abs{y}^a)$. This kind of problems and general properties of elliptic equations with $A_2$ weights are studied in \cite{FKS}.

We now address the nondivergence case.

\begin{lemma} \label{lem:guti}
Given a continuous function $g$ on $\bdary B_R$, such that $g(x,z) = g(x,-z)$, there exists a unique function $u \in C(\overline B_R)$ such that:
\begin{enumerate}
\item[i.] $u$ solves $\partxx u + \abs{z}^\alpha \partzz u = 0$ in $B_R \cup \{z \neq 0\}$ in the classical sense.
\item[ii.] $u \in C^1(B_R)$
\item[iii.] $\partz u(x,0) = 0$
\end{enumerate}
Moreover, for this solution the Harnack inequality result of \cite{CG} applies.
\end{lemma}

\begin{proof}
We point out that any viscosity solution of \eqref{eq:withz} is $C^2$ away from $z=0$, so it would be a solution in the classical sense in $B_R \cup \{ z \neq 0 \}$.

Let us prove uniqueness first. Let us suppose there were two solutions $u$ and $v$ satisfying i., ii., and iii.. For an arbitrary $\eps>0$, consider $w = u - v + \eps \abs{z}$. Then $w \leq \eps R$ on $\bdary B_R$. Let us suppose that $w$ has an interior maximum at a point $x \in B_R$. This point cannot be in $B_R \cap \{z \neq 0\}$, since there $w$ solves the non-degenerate elliptic equation \eqref{eq:withz} and thus it does not have an interior maximum. Therefore, if there is any interior maximum, it has to be on $\{z = 0\} \cap B_R$. But clearly there cannot be a maximum since $\partial_z^+ w > \partial_z^- w$. Therefore $w < \eps R$ in the whole ball $B_R$. Since $\eps$ is arbitrary, we conclude that $u \leq v$ in $B_R$. Similarly we can obtain that $v \leq u$ in $B_R$, so they must coincide. Notice that we only used (i.) and (ii.) for uniqueness.

Now, let us prove existence. The subtle point here is to show that there is a $C^1$ solution. What we do is to prove a uniform $C^{1,\alpha}$ estimate for solutions to the problem
\begin{align}
u^\eps &= g && \text{on } \bdary B_R \\
\lap_x u^\eps + (\abs{z}+\eps)^\alpha \partzz u^\eps &= 0 && \text{in } \bdary B_R \label{aux:approx}
\end{align}
For any $\eps>0$, by the Schauder theory, this problem has classical solutions. If we have a $C^{1,\eta}$ estimate uniform in $\eps$, we take limit as $\eps \to 0$ and obtain the desired solution.

The solutions $u^\eps$ are uniformly bounded in $L^\infty$ due to the maximum principle for uniformilly elliptic equations. The equation \eqref{aux:approx} has constant coefficients with respect to $x$. This means that we can differentiate the equation in a standard way to obtain uniform estimates for derivatives of any order with respect to $x$ (In terms of fractional Laplacian, this corresponds to the fact that functions such that $(-\lap)^s u =0$ are $C^\infty$). Therefore we have that $\lap_x u^\eps$ is bounded independently of $\eps$ in any smaller ball $B_{(1-\delta/2)R}$.

Since $u^\eps \in C^2(B_R)$ and $u$ is symmetric respect to the hyperplane $z=0$, then $\partz u^\eps(x,0) = 0$.

From the equation \eqref{aux:approx}:
\[ 
\begin{split}
\abs{\partzz u^\eps} &= \frac{\abs{\lap_x u^\eps}}{ (\abs{z}+\eps)^\alpha } \\
& \leq \frac{C}{\abs{z}^\alpha}
\end{split} 
\]

Recall that $\alpha = \frac{-2a}{1-a}$ and since $a \in (-1,1)$, then $\alpha < 1$. We can then integrate $\partzz u^\eps$ for any $x,z$ such that $\abs{x} < 1-2\delta$ and $0<z<1-\delta$.

\[ 
\begin{split}
 \abs{\partz u^\eps (x,z)} &= \abs{\int_0^z \partzz u^\eps (x,s) \dd s } \\ 
 &\leq \int_0^z \frac{C}{\abs{z}^\alpha} \dd s = C z^{1-\alpha}
\end{split}
\]

This shows $u^\eps_y$ is $C^\eta$ in $B_{(1-\delta)R}$ for $\eta= \min(1,1-\alpha)$, independently of $\eps$ for any $\delta$.

So, we can take $\eps \to 0$ and extract a subsequence that converges to the desired solution $u$.

Moreover, for any $\eps$, $u^\eps$ is smooth and satisfies the equation \eqref{aux:approx} that has smooth coefficients and for which the Harnack inequality from \cite{CG} applies. Therefore, the same estimate passes to the limit as $\eps \to 0$ and the solution $u$ of the original problem satisfies Harnack inequality.
\end{proof}

\begin{remark}
We have to be careful if we want to study the \em viscosity \em solutions to \eqref{eq:withz}. The naive definition using $C^2$ test functions would not suffice for uniqueness. If, for example, $\alpha > 0$ and we consider the function $u(x,z) =\abs{z}$, then for any $C^2$ function $\phi$ touching $u$ from below at a point $(x,0)$, 
$ \lap_x \phi(x,0) + 0^\alpha \phi_{zz} (x,0) = \lap_x \phi \leq 0 $ and thus $u(x,z)=\abs{z}$ would be a (non-differentiable) \em viscosity \em solution. However if we test against less regular functions, like $\phi(x,z) = z^{2-\alpha}$, then $u$ would fail to be a supersolution. In the notation of \cite{CCKS}, this corresponds to the distinction between $C$-viscosity solutions and $L^p$-viscosity solutions.
\end{remark}

\section{Harnack and boundary Harnack type estimates}
As an application on how to apply the equations \eqref{eq:withy} or \eqref{eq:withz} to the study of fractional harmonic functions, we prove a Harnack inequality and a boundary Harnack inequality using (local) pde methods.

Harnack inequality is not a new result for the fractional Laplacian. It can actually be proved using direct classical potential methods like in \cite{L}. A boundary Harnack estimate for the fractional Laplacian was first proved in \cite{B} using potential methods.

In this paper we derive the Harnack and boundary Harnack inequality for the fractional Laplacian from the Harnack inequality for singular elliptic equations, either with $A_2$ weights (see \cite{FKS}, \cite{FKJ} and also \cite{Penny}) or for certain classes of nondivergence problems (see \cite{CG}).

\begin{thm} [Harnack inequality] \label{thm:harnack}
Let $f :\R^n \to \R$ be nonnegative such that $(-\lap)^s f = 0$ in $B_r$. Then there is a constant $C$ (depending only on $s$ and dimension) such that
\[ \sup_{B_{r/2}} f \leq C \inf_{B_{r/2}} f \]
\end{thm}

\begin{proof}[Proof using \eqref{eq:withy}]
We consider the extension $u$ of $f$ that solves \eqref{eq:withy}. This function $u$ is going to be nonnegative since $f$ is. We reflect it through the $\{y=0\}$ hyperplane. Since $(-\lap)^s f = 0$ in $B_r$, by Lemma \ref{lem:dref}, $u$ is a solution to
\[ \dv \left( |y|^a \grad u \right) = 0 \]
in the $n+1$ dimensional ball of radius $r$ centered at the origin.

We can then apply the result of \cite{FKS} to obtain Harnack inequality for $u$ and thus also for $f$ in half of the ball.
\end{proof}

\begin{proof}[Proof using \eqref{eq:withz}]
We consider the extension $u$ of $f$ that solves \eqref{eq:withz}. This function $u$ is going to be nonnegative since $f$ is. We reflect it through the $\{z=0\}$ hyperplane. Since $(-\lap)^s f = 0$ in $B_r$, then $u$ satisfies the conditions of Lemma \ref{lem:guti}. We can then apply the result of \cite{CG} to obtain Harnack inequality for $u$ and thus also for $f$ in half of the ball.
\end{proof}

\begin{remark}
In the two proofs above we can observe that what is needed for Harnack inequality is that the function $u$ is nonnegative in an $n+1$ dimensional ball. The condition $f \geq 0$ is indeed sufficient for that, but it is not strictly necessary.
\end{remark}

\begin{thm} [Boundary Harnack] \label{thm:bharnack}
Let $f,g : \R^n \to \R$ be two nonnegative functions such that $(-\lap)^s f = (-\lap)^s g = 0$ in a domain $\Omega$. Suppose that $x_0 \in \bdary \Omega$, $f(x) = g(x)=0$ for any $x \in B_1 \setminus \Omega$, and $\bdary \Omega \cap B_1$ is a Lipschitz graph in the direction of $x_1$ with Lipschitz constant less than $1$. Then there is a constant $C$ depending only on dimension such that
\begin{equation} \label{eq:bharnack}
 \sup_{x \in \Omega \cap B_{\frac 1 2}} \frac{f(x)}{g(x)} \leq C \inf_{x \in \Omega \cap B_{\frac 1 2}} \frac{f(x)}{g(x)}
\end{equation}
for any $x,y \in \Omega \cap B_{1/2}(x_0)$.
\end{thm}

\begin{proof}
We consider the extension $u^{(1)}$ of $f$ and $u^{(2)}$ of $g$ that solve \eqref{eq:withy}. As before, we reflect $u^{(k)}$ through $\{y=0\}$ for $k=1,2$ and obtain a solution across this hyperplane through $\Omega$. What we want to do is to find a map that straightens up the domain in order to apply the result of \cite{FKJ}.

First, since $\Omega$ is a Lipschitz domain, we can find a bilipschitz map $\varphi_1: \R^n \to \R^n$ such that $\varphi_1(x_0)=0$ and $\varphi_1(\Omega) \cap B_{1/2} = B_{1/2} \cap \{ x_1 > 0\}$. We can extend this map to $\R^{n+1}$ as constant in the variable $y$. Now, the functions $u^{(k)}_2 = u^{(k)} \circ \varphi_1^{-1}$ are solutions of the equation
\[ \dv( |y|^a b_{ij} \partial_j u^{(k)}_2) = 0 \qquad \text{in $\{y \neq 0\}$ and across $\{y=0\}$ on $\varphi_1(\Omega)$}\]
where the matrix $b_{ij}$ is given by $\frac{D\varphi_1^t D\varphi_1}{\det D\varphi_1}$. Since $\varphi_1$ is bilipschitz, then the singular values of $D\varphi_1$ are bounded below and above and $b_{ij}$ is uniformilly elliptic.

Now we take the map $\varphi_2$ that maps $\R^{n+1} \setminus \{x_1 \leq 0 \wedge y=0\}$ into the half space $\R^{n+1} \cap \{x_1 >0\}$. This map is constant in the variables $x_1,\dots,x_n$, and if we write the pair $(x_1,y)$ in polar coordinates $(r,\theta)$, it maps it to $(r,\theta/2)$. The singular values of $D \varphi_2$ are all equal to $1$ but the one in the direction of $\frac{\partial}{\partial \theta}$ ( i.e. the one in the direction $(-y,0,\dots,0,x_1)$), that is $1/2$. Thus the singular values of $D\varphi_2$ are bounded above and below. The functions $u^{(k)}_3 = u^{(k)}_2 \circ \varphi_2^{-1}$ satisfy the equation
\[ \dv( c_{ij} \partial_j u^{(k)}_3) = 0 \qquad \text{in $B_{1/2} \cap \{x_1 > 0\}$} \]
The matrix $c_{ij}$ is given by $(c_{ij}) = D\varphi_2^t \; (b_{ij}) \; D\varphi_2 \; m(X)$, where the weighted measure $m(X) \dd X$ is the pull back of the measure $|y|^a \dd X$ by $\varphi_2^{-1}$. Since the $A_2$ class of weights is invariant by (local) bilipschitz transformations, the equation still fits into the theory in \cite{FKJ}.

We can apply the boundary Harnack inequality of \cite{FKJ} to the functions $u^{(1)}_3$ and $u^{(2)}_3$ to obtain
\[ \sup_{B_{\frac 1 4} \cap \{x_1 > 0\}} \frac{u^{(1)}_3(x)}{u^{(2)}_3(x)} \leq C \inf_{B_{\frac 1 4} \cap \{x_1 > 0\}} \frac{u^{(1)}_3(x)}{u^{(2)}_3(x)}\]

Recalling that $f(x) = u^{(1)}_3 \circ \varphi_2 \circ \varphi_1 (x_1, \dots, x_n,0)$ and $g(x) = u^{(2)}_3 \circ \varphi_2 \circ \varphi_1 (x_1, \dots, x_n,0)$, we obtain that
\begin{equation} \label{aux:bharnack}
 \sup_{x \in \Omega \cap B_r} \frac{f(x)}{g(x)} \leq C \inf_{x \in \Omega \cap B_r} \frac{f(x)}{g(x)}
\end{equation}
for some universal $r<1$. By a standard covering argument, \eqref{aux:bharnack} together with Theorem \ref{thm:harnack} imply \eqref{eq:bharnack}.
\end{proof}

As before, the condition in Theorem \ref{thm:bharnack} that requires $f$ and $g$ to be nonnegative in the whole $\R^n$ is not sharp. In the proof we can see that what we need to apply the result of \cite{FKJ} is that the corresponding extensions $u^{(1)}$ and $u^{(2)}$ are nonnegative in the unit ball of $n+1$ dimensions. The advantage is that this is a local condition, therefore we can obtain the standard corollary of boundary Harnack saying that $\frac{u^{(1)}}{u^{(2)}}$ is $C^\alpha$ in a neighborhood of the origin. In particular $\frac{f}{g}$ is $C^\alpha$, which is not strictly a corollary of Theorem \ref{thm:bharnack} but of its proof.

\begin{cor}
Let $f$ and $g$ be as in Theorem \ref{thm:bharnack}, then $\frac{f(x)}{g(x)}$ is a $C^\alpha$ function in $\overline{B_{1/2} \cap \Omega}$ for some universal $\alpha$ 
\end{cor}

\begin{remark}
In \cite{FKJ}, the boundary Harnack principle is proven for $A_2$ weights and Lipschitz domains. An extension of that result to nontangentiably accesible (NTA) domains (as in the last section of \cite{ACS}) would lead to a boudary Harnack principle for the fractional Laplacian to domains only requiring a uniform capacity condition.
\end{remark}

\section{Monotonicity Formulas}
Monotonicity formulas are a very powerfull tool in the study of the regularity properties of elliptic PDEs. They have been used in a number of problems to exploit the local properties of the equations by giving information about the blowup configurations. Since the equation \eqref{eq:withy} represents a harmonic function in $n+1+a$ dimensions, it is to be expected that any monotonicity formula known for harmonic functions will have its counterpart for solutions to \eqref{eq:withy}. For example, the simplest one would be that if
\begin{align*}
\Ly u &= 0 && \text{for $(x,y)$} \in B_1^+ \\
\lim_{y \to 0} y^a \party u(x,y)=0 \text{for } |x| \leq 1
\end{align*}
then
\[ \Phi(R) = \frac{1}{R^{n+1+a}} \int_{B_R^+} \abs{\grad u}^2 \ydd X \; \text{  is monotone nondecreasing in $R$}\]
where $B_R^+$ stands for the $n+1$ dimensional half ball $\{(x,y) : \abs{x}^2 + y^2 < 1 \wedge y>0\}$.

The weight $y^a$ is the correct one if we think that the function $u$ is radially symmetric in $a+1$ variables and the $y$ is the modulus of those variables as it is explained in section \ref{sec:rsymmetric}

More interestingly, we have Almgren's frequency formula
\begin{thm} \label{thm:Almgren}
If $u$ is a solution to $\eqref{eq:withy}$ in $B_{R_0}^+$ such that for any $x$ in $B_{R_0}^+$, either $u(x,0)=0$ or $\lim_{y \to 0} y^a \party u(x,y) = 0$, then
\[ \Phi(R) = R  \frac{\int_{B_R^+} \abs{\grad u}^2 \ydd X}{\int_{\bdary B_R^+}  \abs{u}^2 \ydd \sigma} \; \text{  is monotone nondecreasing in $R$ for $R< R_0$}\]

Moreover, $\Phi$ is constant if and only if $u$ is homogeneous.
\end{thm}

The proof is essentially the same as for harmonic functions and it is mainly computational. The original proof can be found in \cite{A}. Here we only need a few minor modifications for our case. First we need the following lemma.

\begin{lemma} \label{lem:Almgren}
Let $u$ be s solution of \eqref{eq:withy} in $B_{R_0} \cap \{y \neq 0\}$, such that $y^a \party u(x,y)$ is bounded and for every $x$ with $\abs{x}<R_0$ either $u(x,0)=0$ or $\lim_{y \to 0} y^a \party u(x,y) = 0$. Then the following identity holds for any $R \leq R_0$.
\begin{equation} \label{eq:diva}
R \int_{\bdary B_R \cap \{y>0\} }  \left( \abs{u_\tau}^2 - \abs{ u_\nu }^2 \right) \ydd \sigma = \int_{B_R^+}  (n+a-1) \abs{\grad u}^2 \ydd x
\end{equation}
where $u_\tau$ stands for the gradient of $u$ tangential to $\bdary B_r$.
\end{lemma}

\begin{proof}
The only thing we have to notice is that since $u$ solves \eqref{eq:withy}, then 
\[ \dv \left( y^a \frac{\abs{ \grad u }^2 }{2} X - y^a \langle X , \grad u \rangle \grad u \right) = y^a \frac{n+a-1}{2} \abs{ \grad u }^2 \]
So we apply divergence theorem in the set $B_R \cap \{y > \eps\}$ to obtain
\[ \begin{split}
R \int_{\bdary B_R \cap \{y> \eps \}} y^a \left( \frac{\abs{\grad u}^2} 2 - \abs{u_\nu}^2 \right) \dd \sigma &+ \int_{\abs{x}<R} -\eps^a \frac{\abs{ \grad u(x,\eps) }^2 }{2} \eps + \eps^a \langle (x,\eps),\grad u(x,\eps) \rangle \party u(x,\eps) \dx \\ &= \int_{B_R^+} \frac{n+a-1}{2} y^a \abs{\grad u}^2 \dx 
\end{split}\]
The second term comes from the integration of the vector field at the bottom of the half sphere. We observe it goes to zero as $\eps \to 0$ since there we have that for each $x$ either $u(x,0)=0$ or $\lim_{y \to 0} y^a \party u(x,y) = 0$ and $a<1$. So we can remove that term and extend the expresion to the whole ball $B_R$.
\[ R \int_{\bdary B_R \cap \{y>0\}}  \left( \frac{\abs{\grad u}^2} 2 - \abs{u_\nu}^2 \right) \ydd \sigma = \int_{B_R^+} \frac{n+a-1}{2}  \abs{\grad u}^2 \ydd x \]
Expanding $\abs{\grad u}^2 = \abs{u_\tau}^2 + \abs{u_\nu}^2$ we obtain \eqref{eq:diva}.
\end{proof}

\begin{proof}[Proof of Theorem \ref{thm:Almgren}]
We will show that $\log \Phi(R)$ is nondecreasing. By the usual scaling argument, it is enough to check it at $r=1$. We compute $\log \Phi$:
\[ \log \Phi(R) = \log R + \log \int_{B_R^+} \abs{\grad u}^2 \ydd X - \log \int_{\bdary B_R^+}  \abs{u}^2 \ydd \theta \]
and now we compute its derivative at $R=1$. Let $S_1 = \bdary B_1 \cap \{y > 0\}$, then
\begin{equation} \label{a1}
(\log \Phi)' (1) = 1 + \frac{ \int_{S_1}  \abs{\grad u}^2 \ydd \sigma} {\int_{B_1^+} \abs{\grad u}^2 \ydd X} - \frac{ \int_{S_1} (2 u u_\nu + (n+a) \abs{u}^2) \ydd \sigma }{\int_{S_1} \abs{u}^2 \ydd \sigma}
\end{equation}

Now we observe that since $\dv ( y^a \grad u)=0$, we have $\dv y^a u \grad u = y^a \abs{\grad u}^2$, and then
\begin{equation} \label{a2}
\begin{split}
 \int_{B_1^+} \abs{\grad u}^2 \ydd X &= \lim_{\eps \to 0^+} \int_{S_1} u u_\nu \ydd \sigma + \int_{\abs{x}<1} u_y u \ydd x \\
 &= \int_{S_1}  u u_\nu \ydd \sigma
\end{split}
\end{equation}

From Lemma \ref{lem:Almgren} we have that
\begin{equation} \label{a3}
\begin{split}
\int_{S_1} \abs{\grad u}^2 \ydd \sigma  &= \int_{S_1 }  \left( \abs{u_\theta}^2 - \abs{ u_\nu }^2 \right) \ydd \sigma + 2\int_{S_1} \abs{u_\nu}^2 \ydd \sigma
\\ &= \int_{B_1^+}  (n+a-1) \abs{\grad u}^2 \ydd X + 2\int_{S_1 } \abs{ u_\nu}^2 \ydd \sigma
\end{split}
\end{equation}

Putting \eqref{a2} and \eqref{a3} together with \eqref{a1} we obtain
\begin{equation*}
\begin{split}
(\log \Phi)' (1) &= 1 + \frac{ \int_{B_1^+} (n+a-1) \abs{\grad u}^2 \ydd X + 2\int_{S_1 }  \abs{ u_\nu}^2 \ydd \sigma} {\int_{S_1}  u u_\nu \ydd \sigma} \\ &\phantom{=} - \frac{ \int_{S_1} (2 u u_\nu + (n+a) \abs{u}^2) \ydd \sigma }{\int_{S_1} \abs{u}^2 \ydd \sigma} \\
&= 1 + (n + a - 1) - (n+a) + \frac{ 2\int_{S_1 } \abs{ u_\nu}^2 \ydd \sigma} {\int_{S_1}  u u_\nu \ydd \sigma} - \frac{ \int_{S_1} 2 u u_\nu \ydd \sigma }{\int_{S_1} \abs{u}^2 \ydd \sigma} \\
&=2 \left( \frac{ \int_{S_1 }  \abs{ u_\nu}^2 \ydd \sigma} {\int_{S_1}  u u_\nu \ydd \sigma} - \frac{ \int_{S_1}  u u_\nu \ydd \sigma }{\int_{S_1}  \abs{u}^2 \ydd \sigma} \right)
\end{split}
\end{equation*}

Therefore $(\log \Phi)' (1) \geq 0$ follows from Cauchy-Schwartz inequality since 
\[ \left( \int_{S_1}  u u_\nu \ydd \sigma \right)^2 \leq 
\left( \int_{S_1}  \abs{u}^2 \ydd \sigma \right) \cdot \left( \int_{S_1 } \abs{ u_\nu}^2 \ydd \sigma \right) \]
and the equality is achieved if and only if $u_\nu = \lambda u$ for some constant $\lambda$. Therefore, only when $u$ is homogeneous of degree $\lambda$.
\end{proof}

\section{Other operators}
The results in this paper suggest that many integro-differential operators can be thought in the same way. Let us suppose we have smooth uniformilly elliptic coefficients $a_{ij}$ in $\R^n \times [0,\infty)$. Consider the following extension problem:
\begin{align}
u(x,0) &= f(x) && \text{for } x \in \R^n \\
\sum a_{ij} u_{ij} (x,y) &= 0 && \text{for } x \in \R^n \text{ and } y > 0
\end{align}
If we now consider the operator $f \mapsto u_y(x,0)$, this is going to be some integro-differential operator of degree one that we can also obtain from a local type pde. We can push the situation even further by considering coefficients $a_{ij}$ that are only measurable, and maybe singular near $y=0$. It is not reasonable to expect that every integro-differential operator can be realized in this way. It looks difficult to characterize the ones that can. But certainly for many it is possible.

Let us consider the case when the operators are invariant under translation in $x$. We want to study the operator that maps the Dirichlet to the Neumann condition for the equation
\[ \lap_x u + a(z) u_{zz} = 0 \]
we take the Fourier transform in $x$ and proceed as in section \ref{sec:ft} to obtain the differential equation
\begin{equation} \label{eq:genwithz}
\begin{aligned}
 -a(z) \phi''(z) + \abs{\xi}^2 \phi(z) &= 0 \\
 \phi(0) &= 1 \\
 \lim_{z \to \infty} \phi(z) &=0
\end{aligned}
\end{equation}
The operator $T : u(x,0) \mapsto -u_z(x,0)$ is then the pseudodifferential operator whose symbol $s(\xi)$ is given by solving \eqref{eq:genwithz} and computing $\phi_z(0)$ for each value of $\xi$.

The question is what symbols $s(\xi)$ can we obtain by the above procedure. We can see that $s(\xi)$ is radially symmetric and monotone increasing in $\abs{\xi}$. Indeed, if $\phi_1$ is a solution of \eqref{eq:genwithz} for $\abs{\xi} = r_0$ and $\varphi_1$ for $\abs{\xi_1} = r_1 \geq r_0$ then 
\[-a(z) \phi_1''(z) + \abs{\xi_0}^2 \phi_1(z) = (\abs{\xi_0}^2 - \abs{\xi_1}^2) \phi_1(z) < 0 \]
Thus $\phi_1 < \phi_0$ by comparison principle.

As $\xi \to 0$, the solution $\phi$ of \eqref{eq:genwithz} converges to $\phi \equiv 1$, so $s(0)=0$. We leave the question if every radially symmetric symbol $s$ such that $s(0)=0$ and it is monotone increasing with respect to $\abs{\xi}$ can be obtained by a suitable choice of $a(z)$.

\section*{Acknowledgements}

Luis Caffarelli was partially supported by NSF grants.

We would like to thank Zoran Vondracek and Renming Song for pointing out a computational mistake in an earlier version of this paper.

\bibliographystyle{plain}   % Here the bibliography 
\bibliography{caffsilv}             % is inserted.
%\index{Bibliography@\emph{Bibliography}}%

\medskip

\noindent Luis Caffarelli \\ 
University of Texas at Austin. \\ 
Austin, TX, USA.\\
\emph{caffarel@math.utexas.edu}

\medskip

\noindent Luis Silvestre \\ 
Courant Institute of Mathematical Sciences. \\ 
New York, NY, USA.\\
\emph{silvestr@cims.nyu.edu}

\end{document}